\documentclass[12pt]{article}
\usepackage{amsmath,amssymb,mathtools,amsfonts}
\usepackage[dvipsnames]{xcolor}
\usepackage{graphicx}
\usepackage{booktabs}
\usepackage{hyperref}
\usepackage[margin=1in]{geometry}

\usepackage{babel}
\usepackage{tikz,pgfplots}
\usetikzlibrary{arrows.meta,calc,decorations.markings,math}
\usetikzlibrary{external}
\usepgfplotslibrary{external}

\usepackage{authblk}

\title{Sparse identification of epidemiological compartment models with conserved quantities}

\author[1*]{Manuchehr Aminian}
\author[2]{Kristin M. Kurianski}
\affil[1]{\small Department of Mathematics \& Statistics, California State Polytechnic University Pomona, 3801 W Temple Ave, Pomona, CA 91768 USA. \url{https://orcid.org/0000-0001-5970-9709}}
\affil[2]{\small Department of Mathematics, California State University Fullerton, 800 N State College Blvd, Fullerton, CA 92831, USA. \url{https://orcid.org/0000-0002-7550-4049}}
\affil[*]{\small Corresponding author. \texttt{maminian@cpp.edu}}

\begin{document}
\maketitle

\begin{abstract}
We investigate the application of a framework for sparse model identification of differential equations from timeseries data in the context of compartmental models in epidemiology. 
Such frameworks often seek a sparse representation from a polynomial basis in the state variables which reproduces the timeseries. 
Out-of-the-box approaches for the underlying sparse regression problem have moderate success reproducing the provided timeseries, but typically fail at producing a consistent, interpretable compartment model and conserving the total population, which are common properties in principled compartment modeling. 
Additionally, the conserved polynomial quantities, such as the sum of state variables, add algebraic nuances to a polynomial design matrix. 
We propose a linear program formulation to solve these issues by posing a pure one-norm objective, sampling from the nullspace of the design matrix, and imposing a set of linear constraints for between-compartment flows.
We conduct several numerical experiments on synthetic data and succeed in ensuring model sparsity, accurately capturing system dynamics, and preserving conservation of population.
\end{abstract}

%


\section{Introduction}\label{sec:intro}

Mathematical epidemiology is usually introduced with the SIR model, tracing back to the work of Kermack and McKendrick \cite{kermack1927contribution}. 
In this work, the authors described a general process with many sequential compartments; for which the three-compartment model of ``susceptible," ``infectious," and ``recovered" individuals is fit to data from the plague epidemic in Mumbai in 1907 \cite{royalsociety1907breport}.
Since then, a veritable zoo of compartment models have been developed and applied to a wide range of infectious diseases, with varying numbers of compartments and  interaction terms to model numerous pathogen and host-specific phenomena at various scales.
Some common extensions modeling a single population are: an inclusion of a loss of immunity (connecting compartment $R$ to $S$); inclusion of an ``exposed" but not infectious compartment; or compartments for vaccination or isolation \cite{hethcote1989three,hethcote2000mathematics, diekmann1990definition, brauer2012mathematical, zaman2008stability,buonomo2008global,pandey2013comparing}.
Including multiple states when exposed, or interactions between multiple species such as through a pathogen's vector (e.g. mosquitoes or fleas) add complexity with the hope of capturing crucial pieces of the dynamical process to incorporate measurable data and make accurate predictions. 
For explaining sub- or super-exponential effects in epidemic dynamics, models also may be extended to graphs (representing restricted travel or friction of travel between well-mixed populations), sometimes revisiting the statistical assumptions of interactions behind infection events \cite{allard2023}, or partial differential equation (PDE) models to capture diffusive or traveling-front effects (see e.g. \cite{vaziry2022}).

Complementary to epidemiology and ecology, there has been a recent surge in the development of tools and investigation of efficacy for purely data-driven methods to discover the underlying modeling relationship or physical law for an observed dynamical system, which in some ways represents an automation of the traditional scientific process.
These can be as simple as discovering the dynamics of a Hookean spring or the classical Lorenz attractor, or as complex as determining underlying PDE models for various phenomena \cite{brunton2016sparse, kaiser2018sparse}. 
Brunton, Proctor, and Kutz (2016) introduced SINDy and serves as a landmark paper to highlight the efficacy of an approach based on an ansatz of an overcomplete dictionary of ``right-hand-side" terms with regression that employs sparsity-promoting regularization \cite{brunton2016sindy}.
Messenger and Bortz (2021) explored the use of a weak-formulation of the differential equations' ansatz to address the vitally important challenge of measurement error in timeseries; where direct methods require derivative estimation and basic methods (such as Forward Euler) are not robust to such noise, weak formulations often are \cite{MESSENGER2021110525}. 
In discovering equations of the forms common in chemical reaction networks and chemical pathways in living systems, Jayadharan et al. (2025) propose a sophisticated method of handling such constraints, where satisfiability of parameter sets lie on general manifolds \cite{jayadharan2025sodas}. 
Our goal in this paper is to specialize tools for data-driven discovery of differential equation models in the context of compartmental models in mathematical epidemiology. 

%
Identifying the number of compartments and their mode of interactions is a persistent challenge in mathematical epidemiology and typically requires subject-matter expertise, brute-force exploration of model and parameter space, or both.
When there is a gap in subject-matter expertise and brute-force is infeasible, there is hope that an automated technique for model discovery informed by data would streamline the model-building process. 
Some recent work in the wake of the COVID-19 pandemic has begun addressing the question of applying the aforementioned sparse model discovery tools to pandemic or epidemic time series data; though their goals and metrics for success vary. 
Horrocks et al. (2020) applied SINDy tools out-of-the-box and found success in qualitative match of various historical epidemics via inclusion of a time-dependent infection parameter in $SI$-style models \cite{horrocks2020algorithmic}.
Jiang et al. (2021) applied their SINDy-LM algorithm to determine model coefficients of a two-dimensional system describing the number of confirmed cases and the number of deaths using COVID-19 data from China, Australia, and Egypt; although modeling challenges related to closed compartment models were outside the scope of their aims \cite{jiang2021modeling}.

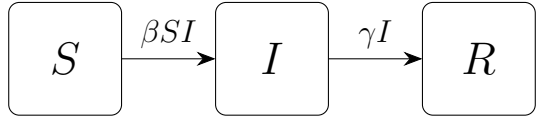
\begin{figure}
\centering
    \tikzsetnextfilename{sir_diagram}
    \begin{tikzpicture}[scale=1] 
        \draw [rounded corners] (1.25,0.5) rectangle (2.75,2); 
        \draw [rounded corners] (4,0.5) rectangle (5.5,2); 
        \draw [rounded corners] (6.75,0.5) rectangle (8.25,2); 
        \draw [-{Stealth[scale=1.5]}] (2.75,1.25) to node [above, pos=0.5] {$\beta SI$} (4,1.25);
        \draw [-{Stealth[scale=1.5]}] (5.5,1.25) to node [above, pos=0.5] {$\gamma I$} (6.75,1.25);
        \node () at (2,1.25) {\Large $S$};
        \node () at (4.75,1.25) {\Large $I$};
        \node () at (7.5,1.25) {\Large $R$};
    \end{tikzpicture}
    \caption{Diagram of SIR compartment model.}
    \label{fig:sir-diagram}
\end{figure}

With model discovery via machine learning (ML) techniques, an ongoing challenge is not necessarily building a model which replicates the original time series, but rather producing models which are explainable and verifiable. 
For example, preserving quantities such as mass, momentum, or energy are important for physics-informed machine learning. 
In the mathematical epidemiology space, model terms represent physical interactions, and conserved quantities relate to dynamics of the overall population. 
In the standard SIR model given by
\begin{equation}
    \begin{array}{l l l l}
         \dot{S}&=&-\beta SI & \\
         \dot{I}&=&+\beta SI & - \gamma I \\
         \dot{R}&=& & +\gamma I,
    \end{array}
    \label{intro-sir}
\end{equation}
the term $\beta SI$ represents the rate of infection, and 
$\gamma I$ is the rate of recovery.
Figure \ref{fig:sir-diagram} shows a diagram of the transitions between compartments.

In the sense of mathematical modeling, these compartment models have two important features. 
First, movement between compartments in this model are symmetric, and the assumption of a closed population yields the identity $\dot{S} + \dot{I} + \dot{R}=0$. 
Most compartment models which exclude long-term population growth have a similar property, or can be easily adjusted to include a compartment tracking population loss to close the system. 
Second, with respect to possible right-hand-side terms, the collection of terms is quite sparse. 
Considering only up to second-degree polynomial terms in the state variables $(S,I,R)$, there are ten possible terms for each equation, thirty in total, of which only four have nonzero coefficients. 
Champion et al. (2020) introduced a sparse identification method that includes constraints; however, the SIR model recovered using their approach still lacks interpretability from a mathematical modeling standpoint \cite{champion2020unified}. 
Since we are interested in model discovery consistent with compartmental modeling principles, we propose applying techniques to modify the optimization problem associated with finding terms in the differential equation with an additional constraint which demands both sparsity and a prescription of a conserved quantity. 

The paper is organized as follows. We introduce our notation and modifications in Section \ref{sec:methods}.
In Section \ref{sec:synthexamples}, we illustrate the procedure in full with examples using synthetic data to demonstrate its efficacy. 
We then conclude in Section \ref{sec:disc} and propose several avenues for ongoing and future work. 

\section{Methods}\label{sec:methods}
In our work, we explore efficacy of constraining the regression problem for discovering the differential equations describing observe dynamics that are compatible with compartmental systems. 

\subsection{Notation and problem formulation} \label{subsec:notation}
Let $x \in \mathbb{R}^M$ with components 
$x=(x_1, x_2, \ldots, x_M)$ be the state vector in autonomous differential equation 
$\dot{x} = f(x)$. We specialize to $f(x)$ which is a linear combination of a dictionary of $\sigma$ functions
$\{\phi_1, \phi_2, \ldots, \phi_\sigma \}$; so for each 
component, $\dot{x_i} = \sum_{j=1}^\sigma W_{ij} \phi_j(x)$ for coefficients
$W_{ij}$. 
The choice of the basis functions and their total number $\sigma$ depends on context; we will typically use a monomial basis of the state variables for polynomials of degree two. 
Though a system of differential equations is usually read in columnar style, 
\begin{equation}
    \begin{array}{c}
         |  \\
         \dot{x_i} \\
         |
    \end{array}
    =
    \begin{array}{c}
    | \\
    f_i(x) \\
    | 
    \end{array},
\end{equation}
in the context of regressing for coefficient matrix $W$, it is more convenient to pose $\Phi(x)=[\phi_1(x), \ldots, \phi_\sigma(x)]$ as a 1-by-$\sigma$ matrix with the coefficient matrix $W \in \mathbb{R}^{\sigma \times M}$. This formulation then poses the system as
\begin{equation}
    \mathrm{-} \dot{x} \mathrm{-} = \Phi(x) W.
\end{equation}
However, we may visualize either $W$ or $W^T$ for convenience.
We call $\Phi$ the design matrix, and 
$\Phi$ will be built using a basis of polynomials for the $M$ state variables $x_i$, 
up to a maximum degree $d$. Therefore the number of functions in the dictionary is 
\begin{equation}
\sigma(M,d) = \sum_{j=0}^{d}\frac{M!}{(M-j)!} = {M+d \choose d}.
\end{equation}
For example, with $M=3$ and $d=2$, the polynomial basis in a lexicographic ordering is
\begin{equation}
\{1, x_1, x_2, x_3, x_1^2, x_1 x_2, x_1x_3, x_2^2, x_2x_3, x_3^2\}.
\end{equation}
Thus, $\Phi(x)$ is a $1\times 10$ matrix of these functions, $W \in \mathbb{R}^{10 \times 3}$ is the coefficient matrix, and $\Phi(x) W$ produces the dynamics in 1-by-3 derivative vector $\dot{x}$.

In the data-driven setting, using a ``strong" formulation, $\Phi(x)$ is replaced by evaluations of $\Phi(x^{(k)})$ for state vector values $k=1,\ldots,N$ observed at $N$ sample points; the derivative $\dot{x}$ is an $N$-by-$M$ matrix and requires an estimation of derivative values; and $W \in \mathbb{R}^{\sigma \times M}$ retains its form. 
We denote the $N$-by-$M$ derivative matrix by $Y\coloneqq\dot{x}$. 
Once $Y$ is formed, the problem can be formulated as a regression problem for $W$: given matrices $Y$ and $\Phi(x)$, find coefficient matrix $W$ which minimizes some combination of normed residual $||Y - \Phi(x)W||$, optimizes coefficient matrix sparsity (e.g. vectorized 1-norm $||W||_1$), and satisfies additional constraint(s) from conserved quantities. 
Since most linear programming expects a formulation in the typical ``$Ax=b$" form with column vectors $x$ and $b$, we vectorize the problem by concatenating columns of $W$. 
Hence, $W \mapsto \widetilde{W} \in \mathbb{R}^{\sigma n \times 1}$, the design matrix is replaced with a block diagonal version with $M$ copies of $\Phi$; i.e. $\Phi \mapsto \widetilde{\Phi} \in \mathbb{R}^{NM \times \sigma M}$, and $Y \mapsto \widetilde{Y} \in \mathbb{R}^{NM \times 1}$, stacked by state variable first.
From this formulation, we investigate constrained approaches to finding approximate solutions to $\widetilde{Y} = \widetilde{\Phi}\widetilde{W}$.

\subsection{Specialization to compartmental models} 
\label{subsec:specialization}
In this subsection we describe an algorithm based on the strong formulation of the problem, based on an ordinary least squares to obtain a candidate solution, followed by seeking a correction to this which both respects the conservation constraint as well as being maximally sparse. 

\subsubsection{Conservation constraints in closed compartment models} \label{subsubsec:conservation}
We consider closed compartment models where
\begin{equation}\label{eq:sumui}
    \sum_{j=1}^M \dot{x_j} = 0.
\end{equation}
Given the matrix form $\dot{x} = \Phi W$, note that 
\begin{equation}
    \dot{x_j} = \sum_{i=1}^\sigma \Phi_{1,i} W_{i,j} \quad\text{ for all }\,j=1,2,\dots, M.
\end{equation}
Thus, a sufficient set of conditions to satisfy \eqref{eq:sumui} is derived by manipulating the matrix form of the expression:
\begin{equation}
    \begin{aligned}
    \sum_{j=1}^M \dot{x_j} &= \sum_{j=1}^M \sum_{i=1}^\sigma \Phi_{1,i} W_{i,j} \\
     0 &= \sum_{j=1}^M \sum_{i=1}^\sigma \Phi_{1,i} W_{i,j} \\
     0 &= \sum_{i=1}^\sigma \Phi_{1,i} \left( \sum_{j=1}^M  W_{i,j} \right) \\
    \Rightarrow 0 &= \sum_{j=1}^M  W_{i,j} \quad \text{for} \quad i=1,\ldots,\sigma;
    \end{aligned}
    \label{W-zero-row-sum}
\end{equation}
where the last equation holds supposing the functions in $\Phi$ form a linearly independent set. 
In other words, each row of $W$ must sum to zero. 
This same calculation holds in the data-driven setting, where derivatives $\dot{x_i}$ are replaced by estimates at each measurement point, and the matrix $\Phi_{1,j}$ in this calculation is replaced by each row of the design matrix; i.e., we will have $N$ copies of the same row-sum constraint \eqref{W-zero-row-sum} on $W$.

In the vectorized version of the system, $\widetilde{W} \in \mathbb{R}^{M \sigma \times 1}$, we have $\sigma$ equality constraints of the form
\begin{equation}
    C_{i,\cdot} \widetilde{W} = 0, \quad \forall i=1,\ldots,\sigma,
\end{equation}
where for fixed $i$, $C_{i,\cdot}$ is a row length $\sigma M$ with ones in entries where $i \equiv j \pmod{\sigma}$ and zeros elsewhere. 
This produces a set of equality constraints on the optimization, represented in matrix equation form,
\begin{equation}
    C \widetilde{W} = 0, \quad C \in \mathbb{R}^{\sigma \times \sigma M}.
    \label{equality-constraint}
\end{equation}

\subsubsection{Selection of a sparse $W$ from candidate} \label{subsubsec:sparse_w}
We propose sparsifying a candidate solution $\widetilde{W}$ which satisfies conservation constraints; $\widetilde{Y}=\widetilde{\Phi} \widetilde{W}$ with $C\widetilde{W}=0$ (for example, linearly-constrained least squares).
We assume in general the design matrix $\Phi$ may have a nontrivial nullspace from which an orthonormal basis can be formed as columns of matrix $\mathcal{N}_{\widetilde{\Phi}}$.
For any $\xi \in \mathrm{col}(\mathcal{N}_{\widetilde{\Phi}})$, we have $\widetilde{\Phi}\xi=0$, therefore $\widetilde{W} + \xi$ also minimizes $||\widetilde{Y} - \widetilde{\Phi} \widetilde{W}||$. 
To aim to increase sparsity in $\widetilde{W}$, we seek a modification $\xi \in \mathrm{span}(\mathcal{N}_{\widetilde{\Phi}})$ to solve the constrained $L_1$ minimization problem:
\begin{equation}
    \begin{aligned}
    \min_{\xi \in  \mathrm{span}(\mathcal{N}_{\widetilde{\Phi}})} || \widetilde{W} + \xi||_1 \\
    \text{subject to } C\xi=-C\widetilde{W}=0.
    \end{aligned}
\end{equation}
To put this in a convenient form for linear programming, let $q = \dim\mathrm{col} (\mathcal{N}_{\widetilde{\Phi}})$ and $z \in \mathbb{R}^q$ be a combination of the columns of $\mathcal{N}_{\widetilde{\Phi}} ( \mathcal{N}_{\widetilde{\Phi}}^T \mathcal{N}_C )$ which samples from constrained solutions, then projects onto the nullspace of the design matrix. 
This bypasses the need for equality constraints downstream. Denote $B\coloneqq \mathcal{N}_{\widetilde{\Phi}} ( \mathcal{N}_{\widetilde{\Phi}}^T \mathcal{N}_C )$. Then the problem reduces to an unconstrained $L_1$ minimization:
\begin{equation}
    \min_{z \in \mathbb{R}^q} || \widetilde{W} + Bz||_1
    \label{L1_min_v1}
\end{equation}
Define slack vector $s \in \mathbb{R}^{M \sigma}$, and for each term $i$ in the objective, 
\begin{equation}
    ||\widetilde{W} + Bz||_1 = \sum_{i=1}^{M \sigma} \big|\widetilde{W} + Bz \big|_i =\sum_{i=1}^{M \sigma} s_i,
    \label{slack-objective}
\end{equation} 
with new objective $c^T \zeta = ( \mathbf{0} \; ; \mathbf{1} )^T (z \; ; s)$, and $2M\sigma$ new inequalities
\begin{equation}
    \begin{aligned}
    \big(\widetilde{W} + Bz \big)_i &\leq s_i, \\
    -\big(\widetilde{W} + Bz \big)_i &\leq s_i.
    \end{aligned}
    \label{slack-inequalities}
\end{equation}
Combining (\ref{L1_min_v1}), (\ref{slack-objective}), and (\ref{slack-inequalities}), this is now a linear program in standard form
\begin{equation}
    \begin{aligned}
     \min_{\zeta \in \mathbb{R}^{q + M\sigma}} \; (\mathbf{0} \; ;  \mathbf{1})^T \, \zeta & \\
    \text{subject to }
    \left( \begin{array}{c | c} B & -I \\ \hline -B & -I \end{array} \right) \zeta 
    &\leq 
    \left( \begin{array}{r} - \widetilde{W} \\ \widetilde{W} \end{array} \right).
    \end{aligned}
    \label{spep-linprog}
\end{equation}
which, after a numerical solution is found, the first block of $\zeta$ is extracted and new solution $\widetilde{W} + Bz$ is formed. 
In the following section we explore the results of this algorithm, specifically, and application of a row-sum constraint, broadly, on synthetic data.

\section{Studies with Synthetic Data}\label{sec:synthexamples}

We now compare the application of out-of-the-box tools to the results of the implementation described in Section \ref{subsec:specialization} with an example SIR model. 
Generically, we generate timeseries $x(t)$ of a known compartment model with coefficient matrix $W$, i.e., $\dot{x} = \Phi(x)W$, then attempt to reconstruct $W$. 
The regressed coefficient matrix is then provided to the same forward simulator so we may consider both timeseries reconstruction (a question of residual norm) and coefficient matrix reconstruction (a question of error norm). 

For the experiments in this section, 
all methods are presented timeseries from the numerical solution of a standard SIR model in Equation (\ref{intro-sir}) with $\beta=0.2$ and $\gamma=0.1$. Solutions are found in Python with an RK45 solver provided with \texttt{scipy.integrate.solve\_ivp}, with relative tolerance $1e-3$ and absolute tolerance $1e-6$. 
Timeseries are sampled $t\in[0,100]$ with a uniform timestep $\Delta t=0.05$. 
When derivatives are needed to build the linear program, either exact derivative information (re-evaluation of $(\dot{S}, \dot{I}, \dot{R})$ at the given state) or forward differencing (FD), i.e. $(S(t_{k+1}) - S(t_k))/\Delta t$ is used. 
The timeseries are provided to one of three workflows: pySINDy with using either an unconstrained or constrained optimizer, scikit-learn's implementation of ElasticNet (an unconstrained regularized regressor), and our methodology presented in \ref{subsec:specialization}. 

\begin{figure}
    \centering
    \includegraphics[width=0.48\linewidth]{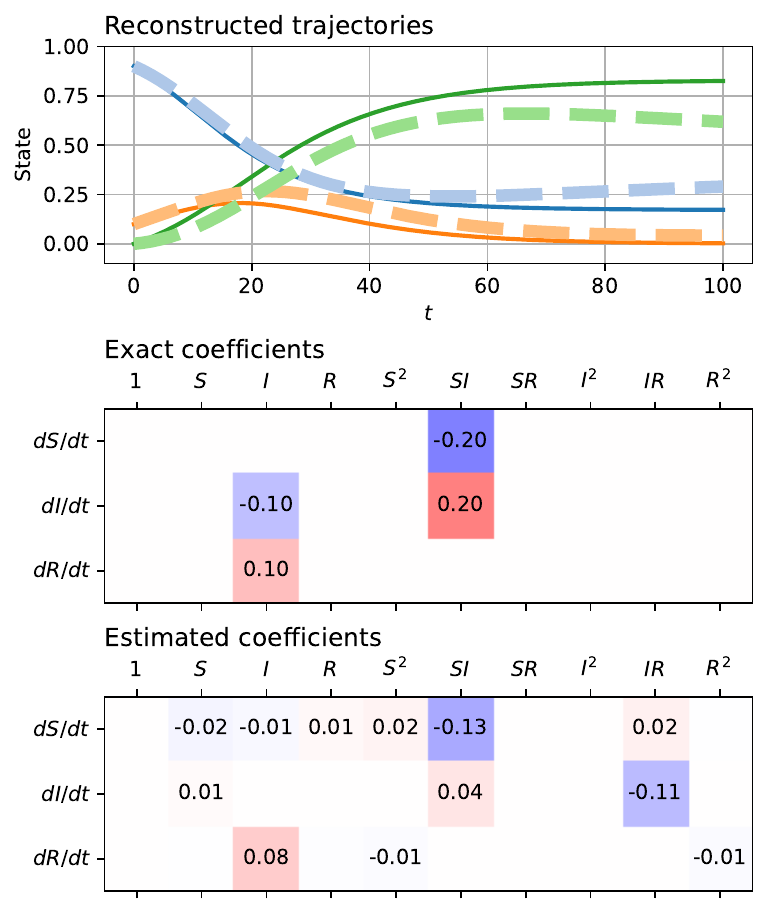}
    \includegraphics[width=0.48\linewidth]{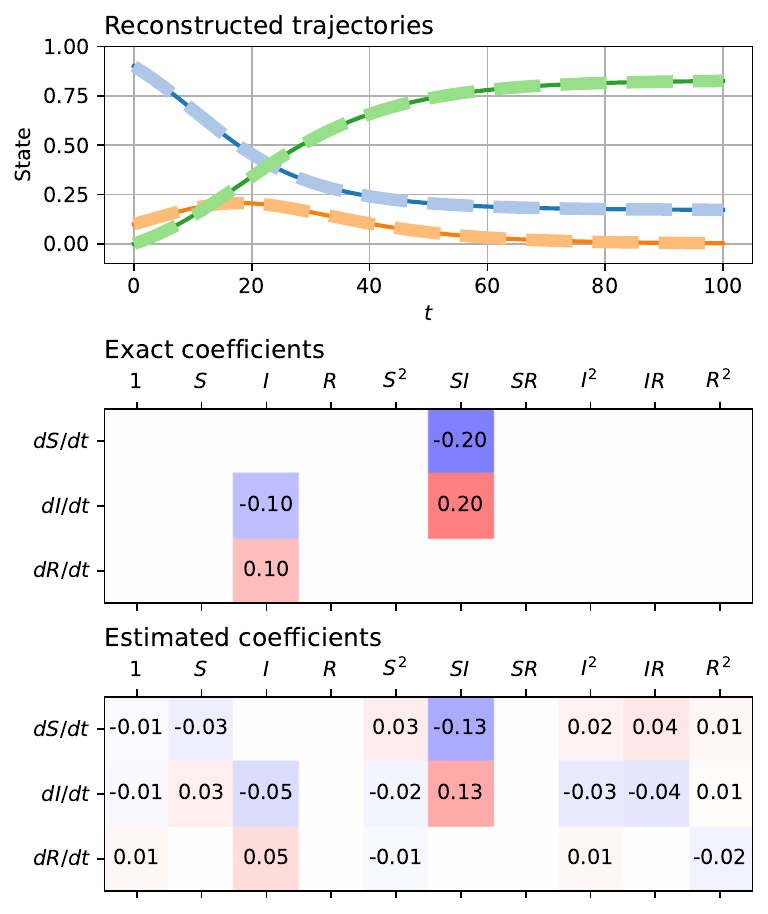}
    \caption{Results of two algorithms, in terms of timeseries reconstruction and recovery of the coefficient matrix. 
    Left column: using forward differences and scikit-learn's implementation of ElasticNet ($\alpha=\mathrm{10^{-6}}$, $\mathrm{l1\_ratio}=0.5$) for regularized sparse regression.
    Right column: applying pySINDy's default settings and the \texttt{ConstrainedSR3} for constrained optimizer with conservation constraints of Section \ref{subsubsec:conservation}. 
    Full results of numerical studies are in Table \ref{tab:numerical_results_exact_derivs}.
    }
    \label{fig:out-of-the-box}
\end{figure}

To begin, Figure \ref{fig:out-of-the-box} illustrates the difficulties out-of-the-box tools have with epidemiological data. 
For visualization purposes, here we display heatmaps of $W^T$ rather than $W$. 
The approach with forward differences and ElasticNet, captures only general qualitative behavior, and fails to match the coefficient matrix (Figure \ref{fig:out-of-the-box}, left).
Implementation with Python package pySINDy with constraint satisfaction via optimizer \texttt{ConstrainedSR3} \cite{champion2020unified} (Figure \ref{fig:out-of-the-box}, right) reproduces the timeseries with good qualitative agreement, though the recovered coefficient matrix $W$ has spurious compartmental interactions, such as an $IR$ term impacting $\dot{S}$ and $\dot{I}$.

\begin{figure}
    \centering
    \includegraphics[width=0.48\linewidth]{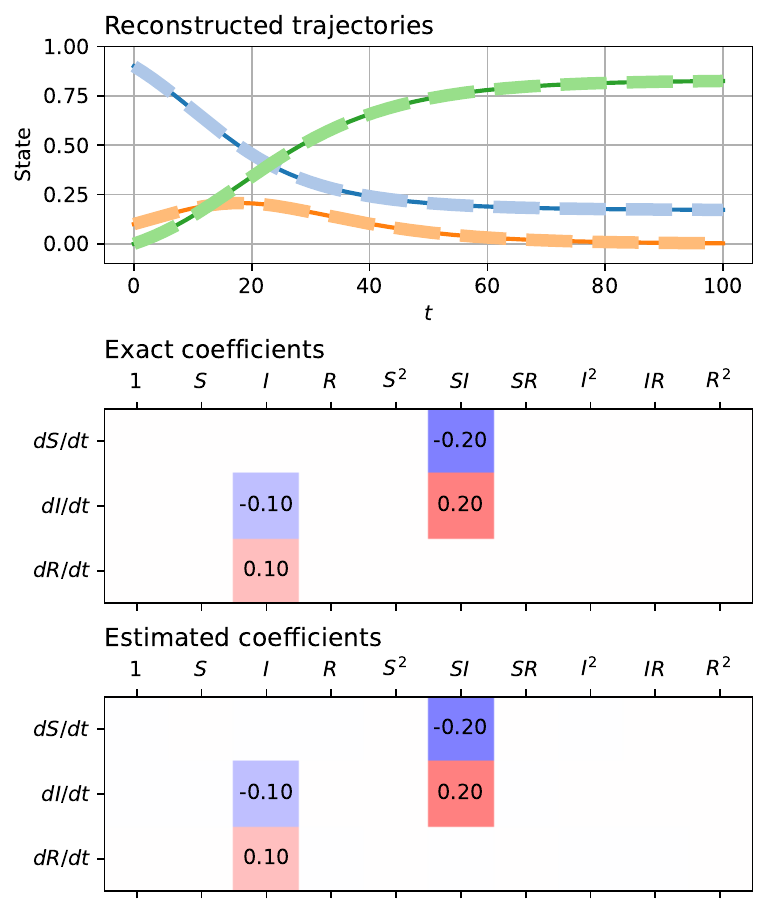}
    \includegraphics[width=0.48\linewidth]{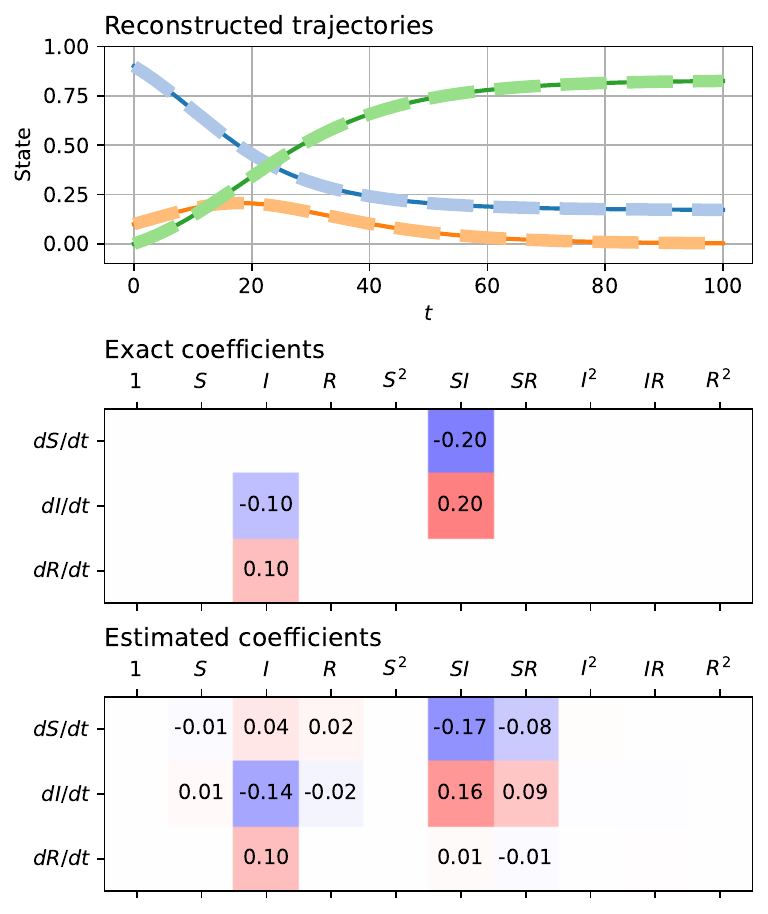}
    \caption{Implementation of algorithm in Section \ref{subsubsec:sparse_w}. 
    Left column: results using exact derivative information recovers both coefficient matrix to at least four digits precision and forward simulation reproduces the timeseries. 
    Right column: when using forward differences to estimate the derivative, quality of the estimated system we recover is degraded. Note coefficients displayed are rounded for readability. 
    }
    \label{fig:ours-synth-sir}
\end{figure}

Figure \ref{fig:ours-synth-sir} demonstrates the advantages of our methodology. 
Using the same synthetic data as in the out-of-the-box implementation, the method outlined in Section \ref{subsec:specialization} produces both a strong qualitative match of the dynamics and quantitative agreement for the coefficient matrix $W$. 
When implementing our methodology, we have a choice in constructing the approximate derivatives $\widetilde{Y}$. 
Figure \ref{fig:ours-synth-sir} (left) shows the results using exact derivative information from the right-hand side of \eqref{intro-sir} to construct $\widetilde{Y}$. 
Figure \ref{fig:ours-synth-sir} (right) shows results using forward differences to construct $\widetilde{Y}$, as was done in Figure \ref{fig:out-of-the-box} (right). 
Either choice of construction for $\widetilde{Y}$ yields qualitative agreement of the timeseries and promising quantitative agreement of the exact and estimated $W$. 
Table \ref{tab:numerical_results_exact_derivs} summarizes the numerical error produced by each algorithm. 
The first three columns evaluate the error in the problem, comparing exact to estimated coefficient matrices $W$. 
The last three columns compare the residuals produced by the timeseries data. 

In all reported measures, the method outlined in \ref{subsec:specialization} using exact derivative information outperforms other algorithms in approximating $W$, most notably ElasticNet and ordinary least squares, which also used exact derivative information. 
In approximating the timeseries, \ref{subsec:specialization} produces smaller errors than ElasticNet and pySINDy, while ordinary least squares produces comparable errors. 
When using forward differences to estimate derivative information, our methodology recovers timeseries information with similar accuracy as pySINDy, while the maximum violation of the conservation constraints from pySINDy are smaller. 
However, the estimated coefficient matrix corresponding to the \ref{subsec:specialization} implementation (Figure \ref{fig:ours-synth-sir}, right) contains fewer spurious compartmental terms than that of the coefficient matrix produced by pySINDy with appropriate constraints (Figure \ref{fig:out-of-the-box}, right).

\begin{table}[h]
    \centering
    \begin{tabular}{l|rrr|rrr}
\toprule
 & mae & 2-norm & maxvio & f\_mae & f\_rmse & f\_mre \\
\hline
ols & 1.1E-01 & 2.2E-01 & 7.4E-16 & 5.0E-14 & 3.1E-14 & 3.4E-13 \\
elasticnet & 1.6E-01 & 2.2E-01 & 9.3E-02 & 1.3E-01 & 7.1E-02 & 2.5E+01 \\
ps\_STLSQ & 1.3E-01 & 2.1E-01 & 4.4E-02 & 4.9E-03 & 2.8E-03 & 9.6E-02 \\
ps\_SR3 & 1.2E-01 & 2.1E-01 & 1.0E-06 & 5.3E-03 & 2.7E-03 & 1.1E-01 \\
\ref{subsec:specialization} exact & 1.7E-14 & 3.4E-14 & 1.3E-14 & 2.2E-13 & 1.4E-13 & 1.4E-12 \\
\ref{subsec:specialization} FD & 9.3E-02 & 1.5E-01 & 1.0E-03 & 4.8E-03 & 2.7E-03 & 1.1E-01 \\
\bottomrule
\end{tabular}

    \caption{Numerical results for several algorithms in reconstructing an SIR system from synthetic data. Algorithms listed include ordinary least squares, scikit-learn's ElasticNet, pySINDy, pySINDy with constrained optimizer, Section \ref{subsec:specialization} with exact derivatives, and Section \ref{subsec:specialization} with derivatives approximated using forward differences. 
    Columns 1-2 give the maximum absolute error and matrix 2-norm of the error for the estimated coefficient matrix $W$. 
    Column 3 gives the maximum violation of the conservation constraint for the row sums of the estimated $W$. Columns 4-6 deal with timeseries reconstruction: the maximum absolute error, root-mean-square error, and the mean relative error respectively.
    All but pySINDy and \ref{subsec:specialization} FD were provided exact derivative information of an SIR time series.
    }
    \label{tab:numerical_results_exact_derivs}
\end{table}

\subsection{Algorithm behavior with a family of SIS models.}
When the design matrix is a polynomial basis of the state variables and the conserved quantity is also a polynomial equation, the data will impose dependencies in the design matrix. 
This fact leads to non-unique solutions to the constrained optimization problem. 
Specifically, the conservation $S + I + R = N$ implies a linear dependence of the design matrix with combination $-N\cdot 1 + 1 \cdot S + 1 \cdot I + 1 \cdot R = 0$. 
This conservation introduces additional algebraic dependencies which leads to further rank-deficiency of the columns of a design matrix of polynomials in $S$, $I$, $R$.
For example, $S(-N + S + I + R)= 0$, and generally, $\mathcal{F}(S,I,R)(-N + S + I + R)=0$ for any polynomial $\mathcal{F}$. In this section, we explore some consequences of these algebraic dependencies via an SIS model. 

\label{subsec:si-exact-minimizer}
\begin{figure}
    \centering
    \includegraphics[width=0.5\linewidth]{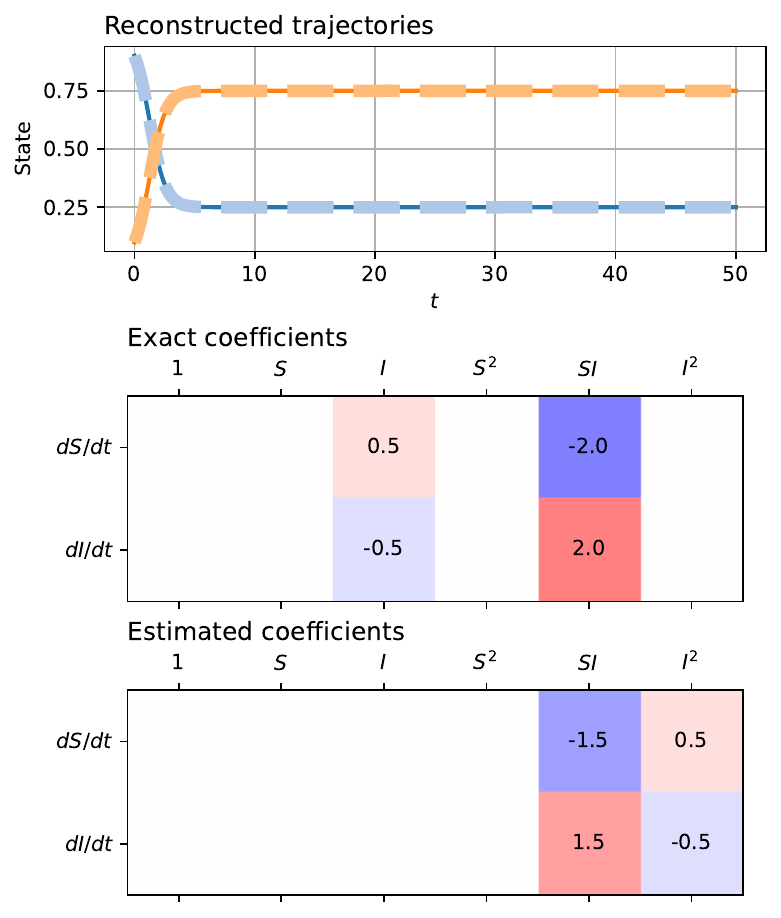}
    \caption{Evaluation of \ref{subsec:specialization} an SIS system. 
    Timeseries are reproduced with a sparse coefficient matrix, but the the systems do not agree termwise. 
    However, the dynamics $(S(t), I(t))$ are provably equivalent by exploiting an assumed conservation equation $S(t)+I(t)=1$.}
    \label{fig:spepi-si}
\end{figure}

To illustrate, we conduct an experiment on an $SIS$ model for which synthetic timeseries data is initialized from the system
\begin{equation}\label{eq:spep-si}
    \begin{aligned}
    \dot{S} &= \frac{1}{2} I - 2SI \\
    \dot{I} &= - \frac{1}{2} I +2SI, \\
    & S(0) + I(0) = 1.
    \end{aligned}
\end{equation}
Figure \ref{fig:spepi-si} illustrates the results of our algorithm on \eqref{eq:spep-si}. 
The dynamics recovered by the algorithm show a different result than the original differential equation: for $S$ this is ${\dot{S} = - \frac{3}{2} SI + \frac{1}{2} I^2}$. 
This is in fact the true minimizer to the problem, which is equivalent to the original differential equation and consistent with the conservation requirement.
To illustrate, since $S(t)+I(t)=1$, substituting $S=1-I$ into \eqref{eq:spep-si} reveals the same dynamics $\dot{I} = 3/2 I -2I^2$. 
A similar analysis can be done by subtracting the two systems and investigating properties of the differences, $\mathcal{S} = \widehat{S} - S$ and $\mathcal{I} = \widehat{I} - I$, where $\widehat{S}$ and $\widehat{I}$ are the reconstructed trajectories produced in Figure \ref{fig:spepi-si}. 
The dynamics for $\mathcal{S}$ are $\dot{\mathcal{S}} = \frac{1}{2} \mathcal{I} ( -1 + \mathcal{S} + \mathcal{I})$, which reduces to $\dot{\mathcal{S}}=0$ for conservative solutions. 
Analytically, we can verify this by seeking a new differential equation of the form
\begin{equation}
\dot{S} = \frac{1}{2} I - 2SI + \mathcal{F}(S,I) (-1+S+I)    
\end{equation} 
since this will produce the same solution paths, assuming $\mathcal{F}$ is polynomial. 
For example, choosing $\mathcal{F} = a + bS + cI$ with undetermined coefficients, we can use a symbolic calculator to find an exact solution for $a$, $b$ and $c$ producing a minimal $\ell_1$-norm of the coefficient matrix $W$ of:
\begin{equation}
    \begin{aligned}
    \dot{S} &= \frac{1}{2} I - 2SI + (a + bS + cI)(-1+S+I) \\
    \dot{I} &= - \frac{1}{2} I +2SI - (a + bS + cI)(-1+S+I)
    \end{aligned}
\end{equation}
%
This is equivalent to minimizing the 1-norm of coefficients of the first equation: 
\begin{equation}
    \min_{a,b,c} |a| + |a-b| + |1/2+a-c| + |b| + |b+c-2|+|b+c|.
\end{equation}
Evaluating this in \texttt{Mathematica 13} produces the unique minimizer $a=b=0$ and $c=1/2$. 
Substitution of these coefficients and collecting terms in the monomial basis produces exactly the numerically discovered system $\dot{S} = -\frac{3}{2} SI + \frac{1}{2}I^2$. 
The 1-norm of the original coefficient matrix is indeed greater ($2(0.5 + 2)=5$) than for the new system ($2(1.5+0.5)=4$), while the number of nonzeros remains the same. 
This confirms the algorithm in Section \ref{subsec:specialization} produced the correct result to the optimization problem, even though the recovered differential equation has purely second-degree terms. 

One can expand this process to a larger class of problems. 
In a similar fashion to the numerical procedure, the 1-norm minimization needed here would itself be converted to a linear program to be solved with symbolic tools, which is typically combinatorially difficult.
On the other hand, this 1-norm minimization is dependent on the original choice of polynomial basis.
Tautologically, one can construct a basis beginning with the correct right-hand-side (e.g. a basis function $\frac{1}{2} I - 2SI$) then expand to complete the function space. This basis will produce a 1-norm minimizer of $(1,0,\ldots,0)$. 
On the other hand, scientific applications often motivate the functional form of terms, such as mass-action dynamics for epidemiological models. 
While these are interesting observations, they are beyond the scope of this work and we postpone them for future work.

\section{Discussion and Future Work}\label{sec:disc}

We present a sparse identification method for discovering differential equation models from epidemiological data that satisfy standard assumptions in compartment modeling. 
Our work addresses two overarching issues: (i) the overabundance of compartment model choices in mathematical epidemiology and (ii) the lack of interpretability from a modeling standpoint of systems produced by many existing sparse identification algorithms. 
We demonstrate the efficacy of our method using synthetic data and compare both accuracy and model interpretability to a host of alternative sparse identification methods. 
While our approach produces similar accuracy to existing methods in terms of timeseries errors, it outperforms current sparse identification models in recovering the compartment model structure of the original system. 
This result holds even when using specialized constrained optimization packages in concert with existing sparse approaches.

Strong formulations for the solution of differential equations, which require direct estimation of derivative values, struggle in the presence of observation error and generally are sensitive to the signal-to-noise ratio (SNR). 
Naive methods such as forward or centered differencing for derivative estimation do not succeed; while pre-processing the timeseries by smoothing may improve the outlook. 
Total variation regularization similarly has been proposed and implemented to preserve the strong formulation while exhibiting robustness to noise with predictable behavior based on the SNR. 
Recent work by \cite{MESSENGER2021110525} among others, identify use of a weak formulation by integrating the ansatz for the differential equation against a set of test functions which may similarly handle noisy data. 
However, the nature of the numerical system and considering dual objectives of timeseries reproduction and conservation satisfaction is sufficiently different compared to the present work.

As seen in Figure \ref{fig:ours-synth-sir}, the coefficient matrix is harder to recover accurately when finite difference estimation is used, even without any direct addition of noise. 
Not shown in this work is some exploration of other higher-precision derivative estimation, which showed no significant improvement on the results. 
We intend to explore more robust derivative estimation in conjunction with these compartment models in the future. 
Another avenue for future studies is the discovery of compartment models from partially observed synthetic data, motivated by the challenges of real-world data collection where more sophisticated compartment models are called for, but include unobservable effects (e.g. susceptible or exposed population).
Finally, a natural extension of our work is applying the methodology to real epidemic data; where both measurement/observation error as well as estimating/inferring the total population $N$ for constructing a closed system is challenging. 
Our ongoing work includes explorations of improved derivative approximations, reconstructing models from unobserved compartments, and applications to real-world data. 

\section{Data and Code}
Code related to this work is available as a public repository at \url{https://github.com/maminian/constrained-epi-model-discovery}. 
This repository includes implementations of the algorithms described in this work in Python, associated tools for analysis, and scripts which generate the figures and tables.


\bibliographystyle{plain}

\end{document}